\documentclass[12pt]{article}
\usepackage{amsmath,amsthm,amssymb,amsfonts}

\usepackage{latexsym}
\usepackage{enumerate}
\usepackage{graphicx}
\usepackage{color}
\usepackage{fullpage}
\usepackage{epsfig}

\newcommand{\NN}{{\mathbb N}}

\newtheorem{definition}{Definition}
\newtheorem{examples}[definition]{Examples}
\newtheorem{theorem}{Theorem}
\newtheorem{lemma}[theorem]{Lemma}
\newtheorem{proposition}[theorem]{Proposition}
\newtheorem{corollary}[theorem]{Corollary}
\newtheorem{questions}{Questions}
\newtheorem{question}[questions]{Question}

\def\pre{{\bf Proof. }}

\begin{document}

\title{On minimal prime graphs and posets}

\author{Maurice Pouzet\footnote{This research, under the auspices of the French-Tunisian
CMCU "outils math\'ematiques pour l'Informatique" 05S1505, was done
while the second author was visiting Sultan Qaboos University. The
support provided by the university is
gratefully acknowledged.}\\
\small ICJ, Universit\'e Claude-Bernard Lyon 1\\
\small 43 Bd. 11 Novembre 1918\\
\small 69622 Villeurbanne cedex France\\
\small and\\
\small University of Calgary \\
\small Department of Mathematics and Statistics\\
\small 2500 University Drive N. W.\\
\small Calgary, Alberta, Canada T2N 1N4\\
\small \and Imed Zaguia\\
\small Department of Mathematics \& Statistics\\
\small Sultan Qaboos University\\
\small P.O.Box 36, Al-Khoud 123. \\
\small Muscat, Sultanate of Oman\\
}

\date{\today}
\maketitle

\footnotetext {\noindent {\it E-mail addresses:}
pouzet@univ-lyon1.fr (M. Pouzet), imed\underline{
}zaguia@hotmail.com (I. Zaguia).}

\begin{abstract}
We show that there are four infinite prime graphs such that every
infinite prime graph with no infinite clique embeds one of these
graphs. We derive a similar result for infinite prime posets with no
infinite chain or no infinite antichain.\\

\medskip \noindent {\bf Keywords:} Prime graph, Prime poset, The neighborhood lattice of a graph,
Incidence structure, Galois lattice, Ramsey Theorem.
\newline {\bf AMS subject classification (2000). 06A06, 06A07}
\end{abstract}

\section{Presentation of the results}\label{sectionone}

This paper is about prime graphs and prime posets. Our notations and
terminology mostly follow \cite {Bo}. The graphs we consider are
undirected, simple and have no loops. That is, a {\it graph} is a
pair $G:=(V, \mathcal E)$, where $\mathcal E$ is a subset of
$[V]^2$, the set of $2$-element subsets of $V$. Elements of $V$ are
the {\it vertices} of $G$ and elements of $\mathcal E$ its {\it
edges}. The {\it complement} of $G$ is the graph ${\overline G}$
whose vertex set is $V$ and edge set ${\overline {\mathcal
E}}:=[V]^2\setminus \mathcal E$. If $A$ is a subset of $V$, the pair
$G_{\restriction A}:=(A, \mathcal E\cap [A]^2)$ is the \emph{graph
induced by $G$ on $A$}. The graph $G$ \emph{embeds} a graph $G'$ and we set $G'\leq G$ if
$G'$ is isomorphic to an induced subgraph of $G$. A subset $A$ of
$V$ is called \emph{autonomous} in $G$ if for every $v \not \in A$,
either $v$ is adjacent to all vertices of $A$ or $v$ is not adjacent
to any vertex of $A$. Clearly, the empty set, the singletons in $V$
and the whole set $V$ are autonomous in $G$; they are called
\emph{trivial}. An undirected graph is called \emph{indecomposable}
if all its autonomous sets are trivial. With this definition, graphs
on a set of size at most two are indecomposable. Also, there are no
indecomposable graph  on a three-element set. An indecomposable
graph with more than three elements will be said \emph{prime}.

The graph $P_4$, the path on four vertices, is prime. In fact, as it
is well known, every prime graph contains an induced $P_4$ (Sumner
\cite{sumner} for finite graphs and Kelly \cite{kelly} for infinite
graphs). Furthermore,  every infinite prime graph contains an
induced countable prime graph \cite{ille}. This leads to the
question: \emph{Which countable prime graphs occur necessarily as
induced subgraphs of infinite prime graphs}?

More specifically, let us say that a graph $G$ is
\emph{minimal prime} if $G$ is prime and every prime induced subgraph with the same
cardinality as $G$ embeds a copy of $G$. One could ask then the
following:

\begin{questions}
\begin{enumerate}[(a)]
\item Does every infinite prime graph embed a countable minimal prime graph?
\item Are there only  finitely many infinite countable minimal prime graphs?
\end{enumerate}
\end{questions}
These questions are the motivation behind this paper. We give a
positive answer for graphs not containing an infinite clique or an
infinite independent set.

\begin{figure}[h]
  \begin{center}
    \includegraphics[width=4.5in]{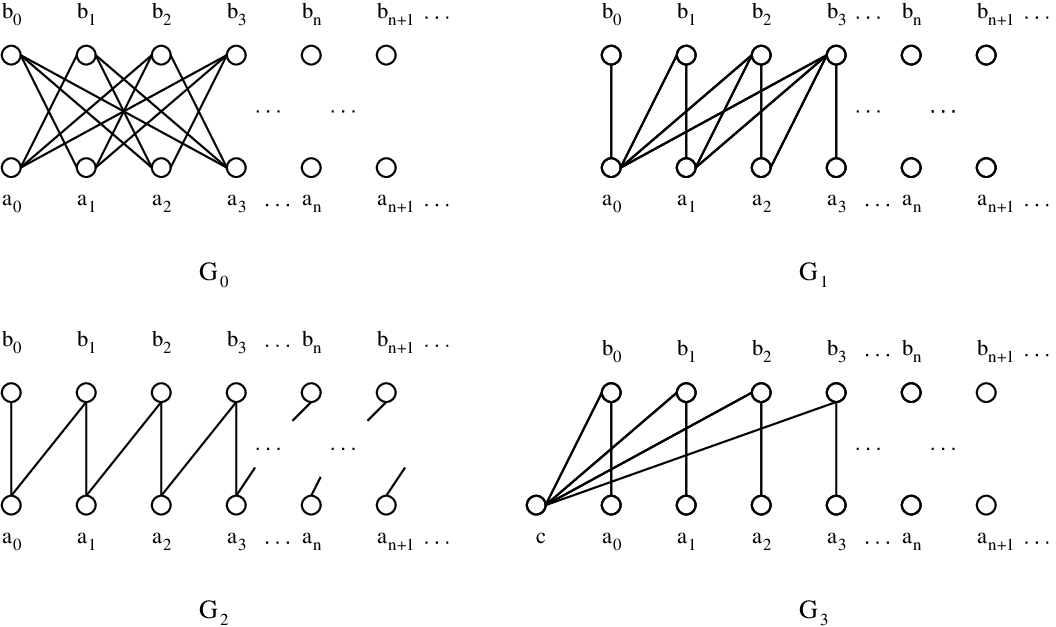}
  \end{center}

  \caption{Minimal prime graphs without an infinite clique.}
  \label{critique}
\end{figure}

In order to state our result, let $\mathcal{G}:=\{G_i : i<4\}$ be the set of graphs defined as
follows.  All these graphs are bipartite, all but $G_3$ have the
same set of vertices which decomposes into two disjoint independent
sets $A:=\{a_i : i\in \NN\}$ and $B:=\{b_i : i\in \NN\}$. A pair
$\{a_i,b_j\}$ is an edge in $G_0$ if $i\neq j$, an  edge in $G_1$ if
$i<j$, an edge in $G_2$ if $j=i$ or $j=i+1$ and, finally,  an edge
in $G_3$ if $j=i$. For $G_3$, a new vertex $c$ adjacent to every
element of $B$ is added to $A\cup B$. The graph $G_0$ is the
comparability graph of $D_{\aleph_{0}}$ (the so-called
\emph{standard poset}, made of the atoms and co-atoms of the Boolean
algebra $\mathcal P(\NN)$ of the subsets of $\NN$), whereas the
graph $G_1$ is the \emph{half complete bipartite graph}. The graph
$G_2$ is the one-way infinite path $P_{\aleph_{0}}$, whereas the
graph $G_3$ is a tree made of countably infinitely many disjoint
edges connected to a single vertex (namely $c$). These graphs are
represented Figure \ref{critique}.

These graphs are prime. A fact which follows from the next
proposition (the proof is easy and let to the reader).

\begin{proposition}\label{primebipartite} A bipartite graph on more than three vertices is prime if and only if
it is connected and distinct vertices have distinct neighborhoods.
\end{proposition}

Moreover, none of these graphs embed in an other. To see that  one
may observe that for each pair $(i,j)$ with $0\leq i\not =j\leq 3$,
there is a finite graph $H_{ij}$ which embeds into $G_i$ and not
$G_j$ (eg take for $H_{01}$ the union of two disjoint edges).


\begin{theorem}\label{main} An infinite prime graph which does not contain an
infinite clique embeds a member of $\mathcal{G}$.
\end{theorem}

An immediate consequence of Theorem \ref{main} (which can be obtained directly)
is that the members of $\mathcal{G}$ are countable minimal prime graphs.

From Theorem \ref{main} we derive two consequences for prime posets.

Throughout, $P :=(V, \leq)$ denotes an ordered set (poset), that is
a set $V$ equipped with a binary relation $\leq$ on $V$ which is
reflexive, antisymmetric and transitive. The \emph{dual} of $P$
denoted $P^{*}$ is the order defined on $V$ as follows: if $x,y\in
V$, then $x\leq y$ in $P^{*}$ if and only if $y\leq x$ in $P$. A
subset $A$ of $V$ is called {\it autonomous} in $P$ if for all
$v\not\in A$ and for all $a,a^{\prime}\in A$

\begin{equation}
(v<a\Rightarrow v < a^{\prime})\;\mathrm{and}\;(a<v\Rightarrow
a^{\prime} < v).
\end{equation}

As for graphs, the empty set, the singletons and the whole set $V$
are autonomous and are said to be {\it trivial}. A poset is {\it
indecomposable} if all its autonomous sets are trivial, it is \emph{prime}
if it is indecomposable with more than three elements.

The \emph{comparability graph}, respectively the
\emph{incomparability graph}, of $P:=(V,\leq)$ is the undirected
graph, denoted by $Comp(P)$, respectively $Inc(P)$, with vertex set
$V$ and edges the pairs $\{u,v\}$ of comparable distinct vertices
(that is, either $u< v$ or $v<u$) respectively incomparable vertices.
We recall the following result (see \cite{kelly}).

\begin{theorem}\label{kelly} A poset $P$ is prime if and only if
$Comp(P)$ is prime. Moveover, if $Comp(P)$ is prime then it has
exactly two transitive orientations, namely $P$ and $P^*$.
\end{theorem}

From this, we have readily:

\begin{proposition}\label {prop:small}
A poset $P$ is minimal prime if and only if $Comp(P)$ is minimal prime.
\end{proposition}

\pre The fact that $Comp(P)$ is minimal prime whenever $P$ is
minimal prime follows directly from the first part of Theorem
\ref{kelly}.\\
The proof of the converse requires also the second part. Indeed, let
$P:=(V,\leq)$ such that $Comp(P)$ is minimal prime. From the first
part of Theorem \ref{kelly} we deduce that $P$ is prime. Let
$V'\subseteq V$ such that $|V|=|V'|$ and $P_{\restriction V'}$ is
prime. Then again $Comp(P_{\restriction V'})$ is prime. Since
$Comp(P_{\restriction V'})=Comp(P)_{\restriction V'}$ and $Comp(P)$
is minimal prime, there is an embedding $f$ of $Comp(P)$ into
$Comp(P)_{\restriction V'}$. If $f$ is not an embedding of $P$ into
$P_{\restriction V'}$, then from the second part of Theorem
\ref{kelly}, $f$ must be order reversing, hence $g=f\circ f$ is an
order preserving map of $P$ into
$P_{\restriction V'}$. \hfill $\Box$\\

As illustrated in Figure \ref{critique}, the four members of
$\mathcal{G}$ are comparability graphs.  Since from Theorem \ref
{main} they are minimal prime, Proposition \ref{prop:small} asserts
that their orientations are minimal prime, whereas Theorem
\ref{kelly} ensures that each one has exactly  two orientations.
Deciding $a_0<b_0$ is each of these graphs we obtain four posets
$Q_0,Q_1,Q_2,Q_3$. The posets $Q_0$ and ${Q_0}^{*}$ are isomorphic
to the standard poset $D_{\aleph_{0}}$. The posets $Q_1$ and
${Q_1}^{*}$ are interval orders, they do not embed in each other.
The posets $Q_2$ and ${Q_2}^{*}$ are two one-way infinite fences,
they are not isomorphic but they do embed in each other. The posets
$Q_3$ and ${Q_3}^{*}$ do not embed in each other. Hence, no member
of $\mathcal{Q}:=\{Q_0,Q_1, {Q_1}^{*}, Q_2, Q_3, {Q_3}^{*}\}$
embeds in another. From Theorem \ref{main} we obtain immediately the
following.

%
%
%
%

\begin{theorem}\label{mainposet1} Every infinite prime poset with no infinite chain embeds a member
of $\mathcal{Q}$.
\end{theorem}

Theorem \ref{main} applies also to incomparability graphs. Indeed,
since a graph is prime if and only if its complement  is prime, a
poset is prime if and  only if its incomparability graph is prime.
We may note that only $G_1$ and $G_2$ are incomparability graphs of
posets. Indeed, as it is well known the comparability graph of a
poset is an incomparability graph if and only if the poset has
dimension at most two \cite{dm}. Since $G_0=Comp(D_{\aleph_{0}})$
and $D_{\aleph_{0}}$ has infinite dimension and since
$G_{3}=Comp(Q_3)$ and $Q_3$  has dimension 3, neither $G_0$ nor
$G_3$ are incomparability graphs. We recall that if a poset $P$ has
dimension 2,  an \emph{order  complement} of $P$ is a transitive
orientation of its incomparability graph. Note that from Theorem
\ref{kelly},  a prime poset has two order complements. The
complements $\overline{G_{1}}$ and $\overline{G_{2}}$ of $G_1$  and
$G_{2}$ have two orientations which are respectively the order
complements $P_1$ and ${P_1}^*$ of $Q_1$ (as well as ${Q_1}^{*}$),
respectively the order complements $P_2$ and ${P_2}^*$ of $Q_2$ (as
well as ${Q_2}^{*}$). See $P_1$ and $P_2$ in Figure \ref
{critique3}. Let $\mathcal{L}$ be the set of these four posets.
These posets are minimal  prime  and  none embeds in an other. From
Theorem \ref{main} we obtain:


\begin{figure}[h]
  \begin{center}
    \includegraphics[width=2.5in]{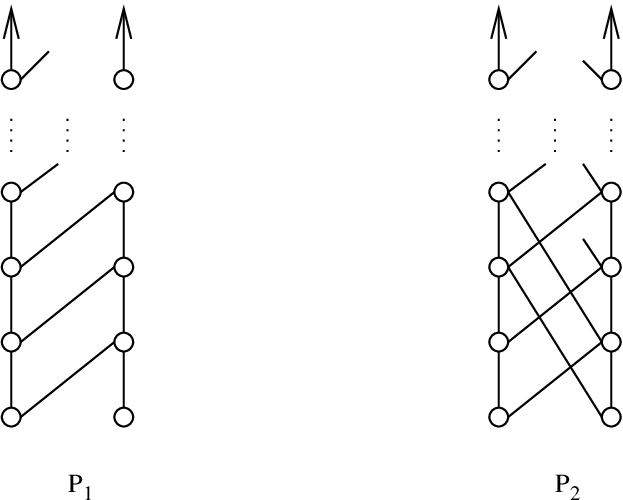}
  \end{center}

  \caption{Minimal prime posets of width two.}
  \label{critique3}
\end{figure}

\begin{theorem}\label{mainposet2} Every infinite prime poset with no infinite antichain embeds a member
of $\mathcal{L}$.
\end{theorem}

Theorem \ref{main} is  a consequence of two properties of the neighborhood lattice of a graph.
To a graph $G:= (V, \mathcal E)$ we associate a complete lattice
$\widehat{N(G)}$, the \emph{neighborhood lattice} of $G$. It is made
of intersections of subsets of $N(G):=\{N_{G}(x) : x\in V
\}$
where $N_{G}(x):=\{y : \{x,y\}\in \mathcal E\}$. Thus ordered by inclusion
this is a complete lattice. As we will see, if $G$ is prime, then
$\widehat{N(G)}$ is infinite. Under the condition that $G$ contains
no infinite clique, we prove that if $\widehat{N(G)}$ contains an
infinite chain, then $G$ embeds $G_0$ or $G_1$. On the other hand if
$\widehat{N(G)}$ contains no infinite chain we prove that $G$ embeds
$G_2$ or $G_3$. Precisely, we prove:

\begin{theorem}\label{maingraph2} Let $G$ be a graph with no infinite clique.
Then $\widehat{N(G)}$ contains an infinite chain if and only if $G$
contains an  induced subgraph isomorphic to $G_0$ or to $G_1$.
\end{theorem}

\begin{theorem}\label{maingraph3} Let $G$ be an infinite prime graph. If all chains in $\widehat{N(G)}$
are finite, then $G$ embeds $G_2$ or $G_3$.
\end{theorem}
The proofs of Theorem \ref{maingraph2} and  Theorem \ref{maingraph3}
rely on properties of incidence structures and on Ramsey's theorem,
with a technique which appeared in \cite {DPR} and \cite{psz}. They
are given in Section \ref{section:maingraph2} and Section
\ref{section:maingraph3}. The properties  we need in order to prove
Theorems \ref{maingraph2} and \ref{maingraph3} are given in the next
section.

We should note that there are other minimal prime graphs and posets. Up to now, we have shown that:

\begin{theorem}\label{thm:examples}There are at least sixteen, respectively twenty two,
countable minimal prime graphs, resp. posets, none embedding in an
other. Furthermore, for every uncountable cardinal $\kappa$, there
are at least fourteen, respectively nineteen, minimal prime graphs,
respectively posets of size $\kappa$, none embedding in an other.
\end{theorem}

The examples leading to Theorem \ref{thm:examples} and the proof are
presented in Section \ref{sample}. In the last section, we present
some questions.

\section{The neighborhood lattice}

Properties of the neighborhood lattice are better understood in
terms of incidence structure and Galois lattices. In the following
subsection, we recall some fundamental properties of these objets.

\subsection{Incidence structures, Galois lattices and coding}

Let $E,F$ be two sets. A {\it binary relation} from
$E$ to $F$ is any subset $\rho$ of the cartesian
product $E\times F$. As usual, we denote by
$x\rho y$ the fact that
$(x,y)\in \rho$ and by $x \neg \rho y$ the negation. The triple
$R:=(E,\rho ,F)$ is an {\it incidence structure}; its
 {\it complement} is $\neg R
:=(E, \neg \rho,F)$, where $\neg \rho:= (E\times F)\setminus \rho$,
whereas its  {\it dual } is $R^{-1} :=(F,\rho ^{-1},E)$, where $\rho
^{-1} := \{ (y,x) : (x,y) \in \rho \}$. For $x\in E$ we set $R(x):=
\{y\in F: x\rho y\}$. Hence for $y\in F$, $R^{-1}(y)=\{x\in X :
y\rho^{-1} x\}=\{x\in X : x\rho y\}$. Let  $R[E]:= \{R(x): x\in E\}$
and $R^{-1}[F]:=\{R^{-1}(y): y\in F\}$. We  denote by  $Gal(R)$ the
set of all intersections of members of $R^{-1}[F]$ (with the
convention that $E\in Gal(R)$). Ordered by inclusion, $Gal(R)$ is a
complete lattice, called the {\it Galois lattice} of $R$. The Galois
lattice of $R^{-1}$ is the set $Gal(R^{-1})$ of all intersections of
members of $R[E]$, ordered by inclusion. A fundamental result about
Galois lattices is:

\begin{theorem}\label{thmgalois}
$Gal(R^{-1})$ is isomorphic to $Gal(R)^*$, the dual of $Gal(R)$.
\end{theorem}

We recall that an incidence structure $R:= (E,\rho ,F)$ is
\emph{Ferrers} if $x\rho y$ and $x'\rho y'$ imply $x\rho y'$ or
$x'\rho y$ for all $x,x'\in E, y,y'\in F$ \cite{ferrers}.
Equivalently, $Gal(R)$ is a chain. We also recall  that a poset $P$
is an interval order iff $(P,<,P)$ is Ferrers.

Let $R:=(E,\rho ,F)$, $R':=(E',\rho ',F')$ be two incidence
structures, a {\it coding from} $R$ {\it into} $R'$ is a pair of
maps $f: E\rightarrow E' , \ \ g: F\rightarrow F'$ such that
$$x\rho y \Longleftrightarrow f(x)\rho ' g(y).$$
When such a pair exists, we say that $R$ has  a {\it coding into}
$R'$.

Bouchet's Coding theorem (\cite {bouchetetat}, see also \cite
{bouchet}) relates the notions of coding and embedding. A
straightforward consequence is this.

\begin {lemma} \label{bouchet2}
  If an incidence structure $R$ has a coding into $R'$, then  $Gal(R)$ embeds into $Gal(R')$.
\end{lemma}

For an example, $Gal((\NN ,\neq ,\NN))= \mathcal P(\NN)$, whereas
$Gal((\NN, <, \NN))$ is the set $I(\NN)$ of initial segments of $\NN$ ordered
by inclusion and $Gal((\NN, >, \NN))$ is the set of final segments
of $\NN$ ordered by inclusion. Hence, from Lemma \ref{bouchet2}, if
one of these structures has a coding in an incidence structure $R$,
the Galois lattice of $R$ embeds the corresponding Galois lattice,
thus contains an infinite antichain. The converse was proved in
\cite{psz}(see Theorem 2.9).

\begin{theorem}\label {thm:psz} The Galois lattice $Gal(R)$  contains an infinite chain if and only if there is a
coding of one of the following incidence structures: $(\NN,=,\NN)$, $(\NN, <, \NN)$ or $(\NN, >, \NN)$ into $R$.
\end{theorem}

Our proof of Theorem \ref{maingraph2} follows similar lines. In
fact, \cite {psz} contains part of Theorem \ref{maingraph2} (see
Corollary 2.15).

\subsection{Basic facts about the neighborhood lattice}

Let $G:=(V,\mathcal{E})$ be an undirected graph without loops. If
$x,y\in V$ we denote by $x\sim y$ the fact that $\{x,y\}\in
\mathcal{E}$ and $x\nsim y$ otherwise. We set $N_G(y):=\{x\in V :
x\sim y\}$. This is the \emph{neighborhood} of $y$. We insist on the
fact that $y\not \in N_{G}(y)$. The \emph{degree} of $y$ in $G$ is
$d_G(y):=|N_G(y)|$, the cardinality of $N_G(y)$. We set
$N(G):=\{N_G(y) : y\in V\}$. Let $\widehat{N(G)}$ be the set of
intersections of subsets of $N(G)$. We make the convention that $V$
is the intersection of the empty set, hence $V\in \widehat{N(G)}$.
Since $\widehat{N(G)}$ is closed under intersections, once ordered
by inclusion this is a complete lattice. We call it the
\emph{neighborhood lattice} of $G$.

%

\begin{figure}[h]
  \begin{center}
    \includegraphics[width=6in]{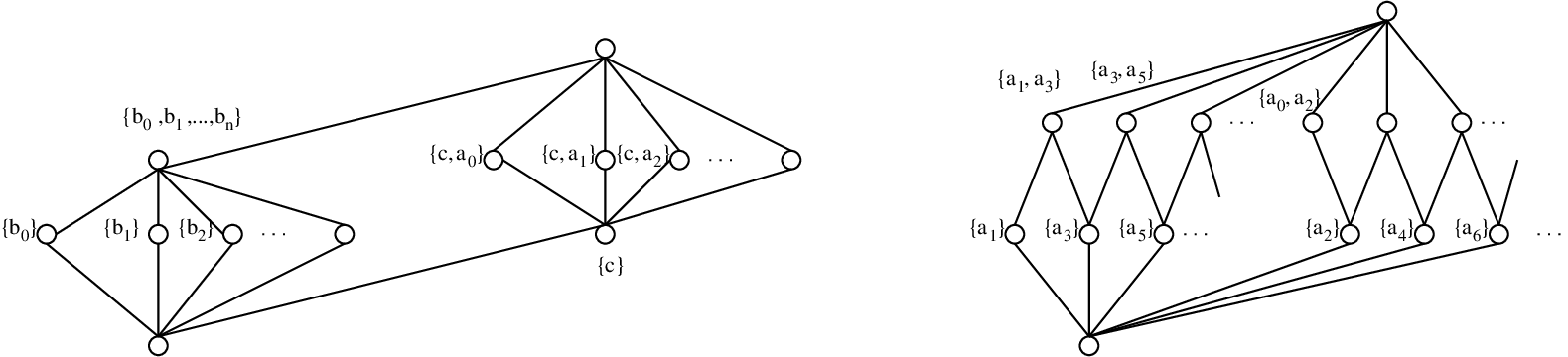}
  \end{center}

  \caption{The neighborhood lattices of $G_2$ and $G_3$.}
  \label{galat}
\end{figure}

Identifying $\mathcal{E}$ to a subset of $V\times V$, or more
precisely setting $\overline{\mathcal{E}}=\{(x,y) : \{x,y\}\in
\mathcal{E}\}$, we have $R^{-1}(y)=N_G(y)$.

\begin{lemma}\label{lem5} The lattice $\widehat{N(G)}$ is the Galois lattice of
$R:=(V,\overline{\mathcal{E}},V)$.
\end{lemma}

Since $R=R^{-1}$, Theorem \ref{thmgalois} yields:

\begin{lemma}\label{lem3} The lattice $\widehat{N(G)}$ is isomorphic
to its dual.
\end{lemma}


Corollary \ref{bouchet2} translates to:

\begin {lemma} \label{bouchet3}
Let $R':=(E',\rho ',F')$ be an incidence structure. If $R'$ has a
coding into $(V,\overline{\mathcal{E}},V)$, where
$G:=(V,\mathcal{E})$ is a graph, then $Gal(R')$ embeds into
$\widehat{N(G)}$.
\end{lemma}
 If  $f$ is an embedding from $G'$ into  $G$ then $(f,f)$ is a coding from $R'$ to $R$. Thus:

  \begin{corollary}\label{cor2} If a graph $G'$ embeds into $G$, then
$\widehat{N(G')}$ embeds into $\widehat{N(G)}$.
\end{corollary}

Lemma \ref{bouchet3} yields:

\begin{corollary}\label{cor3} If a graph $G$ contains an infinite
clique or an induced subgraph isomorphic to $G_0$, then
$\widehat{N(G')}$ contains an induced poset isomorphic to
$\mathcal{P}(\NN)$ ordered by inclusion. If it contains an induced
subgraph isomorphic to $G_1$ then it contains a chain of type
$\omega$ and a chain of type $\omega^*$.\end{corollary}

\pre Let  $G:= (V, \mathcal{E})$. If $G$ contains an infinite clique
or an induced subgraph isomorphic to $G_0$ then  there is a coding
from $(\NN, \not =, \NN)$ in $R:=(V,\overline{\mathcal{E}},V)$,
whereas if  it contains an induced subgraph isomorphic to $G_1$
there is coding from $(\NN,\leq \NN)$ into $R$ and then a coding
from $(\NN, < \NN)$ into $R$. According to Lemma \ref{bouchet3}, in
the first case, $\widehat{N(G)}= Gal(R)$ embeds $Gal((\NN, \not =,
\NN))$ whereas in the second case, $Gal ( G)$ embeds $Gal((\NN, <
\NN))$. Since $Gal((\NN, \not =, \NN))= \mathcal P(\NN)$ and
$Ga((\NN, < \NN))= I(\NN)$, the conclusion  follows. \hfill $\Box$\\

In the sequel, given a subset $X$ of $V$, we set $X^+:=\cap
\{N_G(x): x\in X\}$; eg. $\{x\}^+=N_G(x)$.  With the convention
above, if $X=\emptyset$, then $X^+=V$. Clearly $\widehat{N(G)}=
\{X^+: X\subseteq V\}$. We also set $X^{++}:=(X{^+})^{+}$.

\begin{lemma}\label{lem1} The empty set is the least element of
$\widehat{N(G)}$.

\end{lemma}

\pre $\emptyset=V^+$.\hfill $\Box$\\

A graph $G:=(V, \mathcal E)$ is \emph{point determining} if $x\neq
y$ implies $N_G(x)\neq N_G(y)$ for all $x,y\in V$ (cf.
\cite{sumner2}). We reduce our study to the case of point
determining graphs. Indeed, let $x,y\in V$. Set $x\equiv y$ if
$N_G(x)= N_G(y)$. The relation $\equiv$ is an equivalence relation.
Since $x\not \in N_G(x)$ for all $x$, $N_G(x)=N_G(y)$ implies
$x\nsim y$. Hence, each equivalence class is an independent subset
of $V$. In fact, each equivalence class is also an autonomous subset
in $G$. Set $V/\equiv$ be the set of these equivalence classes and
$p:V\rightarrow V/\equiv$ the map associating to each vertex $x$ its
equivalence class $p(x)$. Let $G/\equiv:=(V/\equiv,E/\equiv)$ where
$E/\equiv:=\{\{p(x),p(y)\} : \{x,y\}\in E\}$. Since the equivalence
classes are independent sets, $G/\equiv$ is an undirected graph with
no loops. Furthermore, $G$ is the lexicographical sum of its
equivalence classes indexed by $G/\equiv$. From this fact follows
readily that $G/\equiv$ is point determining (in fact $\equiv$ is
the unique equivalence relation on $V$ for which the equivalence
classes are independent and autonomous and $G/\equiv$ is point
determining). Furthermore, the map $p$ induces an order isomorphism
from $\widehat{N(G)}$ onto $\widehat{N(G/\equiv)}$. In conclusion,

\begin{lemma}For every graph $G$, $\widehat{N(G)}$ is isomorphic to
$\widehat{N(G')}$ where $G'$ is point determining. In particular,
$\widehat{N(G)}$ and $\widehat{N(G')}$ have the same cardinality.
\end{lemma}

In the remainder of this section, we consider a point determining
graph $G:=(V,E)$.

\begin{lemma}\label{lem2} If $X$ is minimal above $\emptyset$ in
$\widehat{N(G)}$, then $X$ is a singleton.

\end{lemma}

\pre\\ \textbf{Claim 1:} For every $x\in X$ and for every $y\in V$
with
$y\sim x$ we have $X\subseteq N(y)$.\\
Indeed, otherwise set $Y:=X\cap N(y)$. We have $X>Y\cap N(y)\in
\widehat{N(G)}\setminus \{\emptyset\}$ contradicting the minimality
of $X$.\\
From Claim 1 we get:\\
\textbf{Claim 2:} $X$ is an independent and autonomous set.\\
We may now conclude that $X$ is a singleton. Suppose the contrary
and let $x\neq x'$ in $X$. Since $G$ is point determining, $N(x)\neq
N(x')$. Let $y$ be a vertex witnessing this fact. Without loss of
generality we may suppose that $y\sim x$ and $y\nsim x'$. This
contradicts Claim 1.\hfill $\Box$\\

Let $X\in \widehat{N(G)}$, we denote by $\uparrow{X}$ the
\emph{final segment generated by} $X$, that is,
$\uparrow{X}:=\{X'\in \widehat{N(G)} : X\subseteq X' \}$.

\begin{lemma} \label{lem:key}Let $X\in \widehat{N(G)}$ such that
\begin{enumerate}[(1)]
\item $\uparrow{X}$ is infinite.
\item $\uparrow{X'}$ is finite for all $X'\in \widehat{N(G)}$ which contains strictly $X$.
\end{enumerate}

Then
\begin{enumerate}[(1')]
\item $X^+$ is finite.
\item $(X\cup \{x\})^+$ is finite for every $x\not \in X$.
\end{enumerate}

\end{lemma}

\pre Let $x\not \in X$ and set $X':=(X\cup \{x\})^{++}$ and . Since
$X\subseteq X\cup \{x\}$, we have $X\subseteq X'$ and since $x\in
X'\setminus X$, $X\neq X'$. Thus from (2) $\uparrow X'$ is finite.
Since $X'\subseteq \{y\}^+=N_G(y)$ for all $y\in (\{x\}\cup X)^+$,
the set $\{N_G(y) : y\in (\{x\}\cup X)^+\}$ is finite. Since $G$ is
point determining, $(\{x\}\cup X)^+$ is finite as required.\hfill
$\Box$

\begin{corollary}The following properties are equivalent.
\begin{enumerate}[(1)]
\item For every $X\in \widehat{N(G)}\setminus \{\emptyset\}$,
$\uparrow{X}$ is finite,
\item $N_G(x)$ is finite for every $x \in V$.
\end{enumerate}
\end{corollary}

\pre $(1)\Rightarrow (2)$ Apply Lemma \ref{lem:key} with $X:=\emptyset$.\\
$(2)\Rightarrow (1)$ Let $X$ be non empty. The set $\uparrow{X}$ is
finite if there are only finitely many $N_{G}(y)$ containing $X$.
Pick $x\in X$. Since $N_G(x)$ is finite, the numbers of $N_G(y)$
such that $y\in N_G(x)$ is finite. In particular the number of
$N_G(y)$ containing $X$ is finite. \hfill $\Box$

\section{Proof of Theorem \ref{maingraph2}}\label {section:maingraph2}
If $G$ contains an induced subgraph isomorphic to $G_0$ or to $G_1$
then, according to Corollary \ref{cor3}, $\widehat{N(G)}$ contains an
infinite chain. Conversely, suppose that $\widehat{N(G)}$ contains
an infinite chain. Then   it contains a chain of type $\omega$ or
$\omega^*$. Since the lattice $\widehat{N(G)}$ is selfdual (Lemma
\ref {lem3}), it contains a chain of type $\omega$, meaning that
there exists a strictly increasing sequence $(X_n)_{n\geq 0}$ of
members of $\widehat{N(G)}$. From this,  we may define two maps
$f_0:\NN\rightarrow V(G)$ and $f_1:\NN\rightarrow V(G)$ such that
for all $n\in \NN$:
\begin{enumerate}[(1)]
\item $f_0(n) \in X_{n+1}$ and $f_1(n)\in {X_{n}}^{+}$.
\item $f_0(n)\nsim f_1(n)$.\\
Indeed, since $X_{n+1}\nsubseteq X_{n}={X_{n}}^{++}$, there are
$a\in X_{n+1}$ and $b\in {X_n}^+$ such that $a\nsim b$. Set
$f_0(n):=a$ and $f_1(n):=b$.\\
Beyond that, the maps $f_0$ and $f_1$ has the
following properties:
\item $f_0(n)\sim f_1(m)$ for all $n<m$.\\
Indeed $X_{n+1}\subseteq X_{m}\subseteq X_{m+1}$ thus
${X_{m+1}}^{+}\subseteq {X_{m}}^{+}\subseteq {X_{n+1}}^{+}$. Since
$f_1(m)\in {X_{m}}^{+}$ we have $b_{m}\in {X_{n+1}}^{+}$. Since
$f_0(n) \in X_{n+1}$ this yields $f_0(n)\sim f_1(m)$.
\item $f_1(n)\neq f_1(m)$ for all $n<m$.
\item $f_0(n)\neq f_0(m)$ for all $n<m$.
\end{enumerate}

Indeed, from (1) we have $f_0(n)\nsim f_1(m)$, thus (3) holds.
Similarly (1) yields $f_0(m)\nsim f_1(m)$ and (2) yields $f_0(n)\sim
f_1(m)$.
Thus (4) holds. The proof of (5) is similar.\\

Let $[\NN]^2$ be the set of two element subsets of $\NN$ identified
with ordered pairs $(n,m)$ such that  $n<m$. Divide $[\NN]^2$ into
blocks such that two such pairs $u:=(n_0,n_1)$ and $u':=(n'_0,n'_1)$
are in the same block if
\begin{equation}
f_i(n_k) \rho f_j(n_l) \Leftrightarrow f_i(n'_k) \rho f_j(n'_l)
\end{equation}
holds for all $i,j,k,l\in \{0,1\}$ and $\rho \in \{=,\sim\}$. \\As
it is easy to see the number of blocks is finite.  Indeed, it is
bounded by $2^{24}$ (each block can be coded by a relational
structure made of six binary relations on a two element set).
Ramsey's theorem on pairs ensures that there is an infinite subset
$I$ of $\NN$ such that all pairs belong to the same block. Let
$\phi(0)<\phi(1)<...<\phi(n)...$ be an enumeration of $I$ and
$\overline{f_i}:=f_i\circ \phi$ ($i<2$). Then equivalence $(2)$
holds with $f_i$ and $f_j$ replaced by $\overline{f_i}$ and
$\overline{f_j}$, meaning that all pairs of $[\NN]^2$ are in the
same block. Thus, without loss of generality, we may choose $f_0$
and $f_1$ such that equivalence (2) holds. We  say then that the
pair $(f_0, f_1)$ behaves \emph{uniformly} on $\NN$. In this case we
have the following additional properties:

\begin{enumerate}[(6)]
\item $f_0(n)\nsim f_0(m)$
\end{enumerate}
\noindent and
\begin{enumerate}[(7)]
\item $f_1(n)\nsim f_1(m)$
\end{enumerate}

\noindent for all $n<m$.
 Indeed, if $f_0(n_{0})\sim f_0(m_{0})$ for some pair
$(n_0,m_0)$ then since $(f_0,f_1)$ behaves uniformly, $f_0(n'_{0})\sim
f_0(m'_{0})$ holds for all other pairs, and thus $G$ contains an
infinite clique. The proof of (7) is similar.

\begin{enumerate}[(8)]
\item $f_0(n)\neq f_1(m)$ for all $n\neq m$.
\end{enumerate}

Indeed, if $n<m$ this follows from (2). If $f_0(n)=f_1(n)$ for some
$n$, then since the pair $f_0, f_1)$ behaves uniformly,
$f_0(n)=f_1(n)$ for all $n$. But for $n<m$ we get $f_0(n)\sim
f_1(m)=f_0(m)$ contradicting (5). If $f_0(n)=f_1(m)$ for some $m<n$,
then $f_0(n)=f_1(m)$ for all $m<n$. Taking $m<n<n+1$ we get
$f_0(n)=a_{n+1}$ contradicting (4).

So far, the sets  $A':=\{f_0(n) : n<\omega\}$ and $B':=\{f_1(m) :
m<\omega\}$ are
two disjoint independent subsets of $G$ for which (1) and (2) hold. We consider two cases.\\
\textbf{Case 1.} (9) $f_1(n)\sim f_0(m)$ for some $n<m$.\\
Since $(f_0,f_1)$ behaves uniformly, this property holds for all $n<m$.
In this case $G_{\restriction  A'\cup B'}$ is isomorphic to $G_0$.\\
\textbf{Case 2.} (10) $f_1(n)\nsim f_0(m)$ for some $n<m$. Again,
this property holds for all $n<m$. In this case $G_{\restriction
A'\cup B'\setminus \{b_0\}}$ is isomorphic to $G_{1}$, via the map
$\phi$ from $G_1$ to $G$ defined by $\phi(a_n)=f_0(n)$ and
$\phi(b_n)=f_1(n+1)$.\hfill $\Box$

\section{Proof of Theorem \ref{maingraph3}}\label {section:maingraph3}

Let $G$ be an infinite prime graph. Since its is prime, it is point
determining. This allows us to apply Lemma \ref{lem:key}. Let
$\mathcal{X}$ be the set of $X\in \widehat{N(G)}$ such that the
final segment $\uparrow{X}$ of $\widehat{N(G)}$ is infinite. Since
$\widehat{N(G)}$ is infinite, $\emptyset \in \mathcal{X}$. Since
$\widehat{N(G)}$ contains no infinite chain, $\mathcal{X}$ has a
maximal element (this may require the  axiom of dependent choices).
Let $X$ be such an element. Then for each $X'\in \widehat{N(G)}$
containing strictly $X$, the final segment $\uparrow {X'}$ is
finite. According to Lemma \ref{lem:key}, $X^+$ is infinite. Let
$G':=G_{\restriction X^+}$.\\
\textbf{Case 1.} $G'$ contains an infinite connected component.\\
In this case, $G'$ contains an infinite path. Indeed, according to
Lemma \ref{lem:key}, for each $x\in X^+$, the degree of $x$ in $G'$ is finite.\\
\textbf{Case 2.} All connected components of $G'$ are finite.\\
In this case, since $G$ is prime, $G$ is connected hence $G'\neq
G$, that is, $X\neq \emptyset$.\\
\textbf{Claim 1.} For every connected component $C$, but perhaps
one, there are $a_C\in C$ and $b_C \in V\setminus X^+$ such that
$a_C\sim b_C$.\\
\textbf{Proof of Claim 1}. Let $C$ be a connected component of $G'$
Since $C\subseteq X^+$, we have $X\subseteq F(a):=N_{G}(a)\cap
(X\setminus X^+)$. If $X=F(a)$ for every $a$, then $C$ is autonomous
in $G$. Since $G$ is prime, $C$ is a singleton, that is, $C=\{a\}$
for some $a$. There is no other connected component $C'$ reduced to
a singleton, because otherwise the set $\{a,a'\}$, where
$A'=\{a'\}$, is autonomous in $G$. Thus, all connected components
$C$ but one contain an element
$a_C$ such that $X\neq F(a_C)$. Pick $b_C\in F(a_C)\setminus X$. $\blacksquare$.\\
We define inductively  two maps $f_0:\NN\rightarrow V(G)$ and
$f_1:\NN\rightarrow V(G)$. Suppose that
$(f_0(i),f_1(i))_{i<n}$ has been defined. According to Lemma
\ref{lem:key},  $(\cup_{i<n}N_{G}(f_0(i))\cap X^+$ is finite. Pick a
connected component $C$ of $G'$ which does not meet
$\cup_{i<n}N_{G}(f_1(i))$ and set $f_0(n):=a_C$ and $f_1(n):=b_C$. The
sequence $(f_0(n),f_1(n))_{n\in \NN}$ has the following properties:
\begin{enumerate}[(1)]
\item $f_1(n)\neq f_0(m)$ for all $n\neq m$.\\ Indeed, $f_0(m)\in X^+$ and
$f_1(n)\not \in X^+$.
\item $f_1(n)\nsim f_0(m)$ for all $n< m$.\\ Indeed, $f_0(m)$ has been selected
in a connected component which does not meet $N_{G}(f_1(n))$.
\item $f_0(n)\neq f_0(m)$ for all $n\neq m$.\\ Indeed, these elements are
chosen in different connected components of $G'$.
\item $f_1(n)\nsim f_1(m)$ for all $n< m$.\\ Otherwise, we would have
$f_1(n)=f_1(m)$ for some $n<m$. Since $f_1(m)\sim f_0(m)$ we would
have $f_1(n)\sim f_0(m)$, contradicting (2).
\end{enumerate}

Apply Ramsey Theorem as in the proof of Theorem \ref{maingraph2}.
There is an infinite subset $I$ of $\NN$ on which the pair
$(f_0,f_1)$ behaves uniformly. Without loss of generality we may
suppose that $I=\NN$ (otherwise relabel $I$ with the integers). From
the fact that $N(G)$ contains no infinite chain, $G$ contains no
infinite clique. This excludes $f_1(n)\sim f_1(m)$. Again the fact
that $\widehat{N(G)}$ contains no infinite chain excludes
$f_0(n)\sim f_1(m)$ for $n<m$. Let $C:=\{c\}\cup \{f_0(n),f_1(n) :
n\in \NN\}$ where $c\in X$. Then $G_{\restriction  C}$ is isomorphic
to $G_{3}$. \hfill $\Box$

\section{Examples of minimal prime graphs and posets}\label{sample}
All the minimal prime graphs that we have been able to obtain so far
have,  at the exception of $G_{2}$ and its complement, a common
feature, that we present in full generality. As a byproduct, we
obtain examples of minimal prime graphs of arbitrarily cardinality.
Then we identify those which are comparability graphs.

\subsection{Uniform graphs}
 A graph $G:=(V,\mathcal{E})$ is \emph{uniform} (\cite {psz}) if
 $V$ is the disjoint
union of a finite set $K$ and a set of the form $E\times \{0,1\}$.
The set $E$ is equipped with a linear order $\leq$. For two distinct
vertices $u$ and $v$, the fact that they form an edge or not only
depends upon how $x$ and $y$ are related by the order and upon the
values of $i$ and $j$ if $u:=(x,i)$ and $v:=(y,j)$ or upon the value
of $i$ if $u:=(x,i)$ and $v\in K$. Formally, this translates to:
\begin{enumerate}[(1)]
\item $(x_k,i)\rho (x_l,j) \Leftrightarrow (x'_k,i)\rho (x'_l,j)$.
\item  $(x_k,i)\rho y \Leftrightarrow (x'_k,i)\rho y$.
\end{enumerate}
for all $x_0<x_1$, $x'_0<x'_1$ in $E$, $y$ in $F$, $i,j,k,l\in
\{0,1\}$, $\rho \in \{\sim,=\}$.

\begin{examples} Let $C:=(E,\leq)$. Set
$G_{k}(C):=(V_k,\mathcal{E}_k)$ for $k\in \{0,1,3,4\}$ with
$V_k:=K_k\cup E\times \{0,1\}$, $K_k=\emptyset$ in case $k\neq 3$ and
$K_3=\{c\}$. Set

\begin{enumerate}[(1)]
\item $(x,i)\sim_0 (x',i')$ if $i\neq i'$ and
$x\neq x'$.
\item $(x,i)\sim_1 (x',i')$ if $i\neq i'$ and
either $i<i'$ and $x\leq x'$ or $i>i'$ and $x'\leq x$.
\item $(x,i)\sim_3 (x',i')$ if $i\neq i'$ and
$x= x'$, $c\sim (x,i)$ if $i=1$.
\item $(x,i)\sim_4 (x',i')$ if either $i= i'=0$ and $x\neq x'$
or $i\neq i'$ and $x= x'$.

\end{enumerate}

We introduce three more graphs, in a more informal way.

\begin{enumerate}[(1)]
\item The graph $K(C)$ has the same set of vertices as $G_1$ and edge
set $\mathcal{E}_1 $ augmented with the set of unordered pairs of
distinct elements of $E\times \{0\}$.
\item The graph $G_c(C)$ is obtained from $K(C)$ by adding a new
vertex $c$ adjacent to all elements of $E\times \{1\}$.
\item $G_{ab}(C)$ is obtained from $K(C)$ by adding two extra
vertices $a$ and $b$ and an edge between $a$ and all elements of
$\{b\}\cup E\times \{0\}$.

\end{enumerate}

\end{examples}

\subsection{Examples of minimal prime graphs}
 Let us recall that a chain $C$ isomorphic to the chain of
nonnegative integers has order type $\omega$. More generally, if $C$
is well ordered its \emph{order type} is the unique ordinal to which
$C$ is isomorphic. In the remainder of this section we will mostly
consider initial ordinals, eg
$\omega,\omega_1,\omega_2,...,\omega_{\omega},...$. These are
cardinal numbers, the \emph{aleph}'s. With the axiom of choice,
there are no others.

In the definition of $G_i(C)$'s, $i \in \{0,3,4\}$,  the order of
$C$ is irrelevant. Replacing $C$ by $\omega$, we  obtain
$G_i=G_i(\omega)$ for $i \in  \{0,1,3, 4\}$. The order is quite
relevant in the other cases.

\begin{theorem}\label{thm:listminimalprime graph}
For every infinite initial ordinal  $\kappa$, the graphs
$G_i(\kappa)$ for $i<5, i\not =2$,  their complements, and the
graphs  $G_c(C), G_{ab}(C), \overline G_{ab}(C)$ for  $C\in
\{\kappa, \kappa^*\}$ are minimal prime. Furthermore, none of these
fourteen  graphs embeds in an other. \end{theorem}

The fact that those  graphs are minimal prime is an immediate consequence of the following lemma.

\begin{lemma} \label {thm:minichain}Let $C:=(E, \leq)$ be a chain.

\begin{enumerate}[(1)]
\item The graphs $G_i(C)$, $i<5, i\not =2$, $G_c(C)$ and $G_{ab}(C)$
are prime, provided that $\vert E\vert \geq 3$ if $i=0$ and $\vert E\vert \geq 2$ in all other cases.
\item If $\kappa$ is an infinite initial ordinal and $C$ is a chain of order
type $\kappa$ or $\kappa^*$, these graphs are minimal prime.
\end{enumerate}

\end{lemma}

\pre Let $G$ be one of the six graphs listed in Theorem \ref {thm:minichain}.
$(1)$ The graph $G$ is prime.  If $G\in \{G_i(C) : i\in \{0,1,3\}$,  apply Proposition
\ref{primebipartite}.\\
Let $G\in \{G_4(C),G_c(C),G_{ab}(C)\}$. Let  $X$ be an autonomous
subset in $G$ with $|X|\geq 2$. To prove that $G$ is prime we need
to prove that $X=V$. For that we make a repeated use of the
observation that if a vertex $v$  \emph{separates} two vertices $u$
and $u'$ of $X$ (that is, $v\sim u$ and $v\nsim u'$ or $v\nsim u$
and $v\sim u'$), then $v\in X$.

Set $X_i:= X\cap  E\times \{i\}$ for $i\in \{0, 1\}$.

{\bf Claim 1}. Let $i, j\in \{0,1\}$ with $i\not =j$. If  $\vert X_i\vert \geq 2$ then
$ X_j\not =\emptyset$.

Indeed, let $u:=(x,i)$, $u':=(x,i)$ be two distinct elements of
$X_i$. We may suppose that $x<x'$. Let $v:=(x',0)$ if $i=1$, $v:=
(x,1)$ if $i=0$. Then $v$ separates $u$ and $u'$. Thus $v\in X_j$,
proving our claim.

{\bf Claim 2}. If $X_0$ and $X_1$ are non empty then $X=V$.

Note first that there is some $z\in E$ such that $(z,i)\in X$ for
all $i\in \{0,1\}$. Indeed, let $(x, 0), (x',1)\in X$. If $ x=x'$ we
are done. If $x\not =x'$ then $(x,1)$ separates $(x,0)$ and $(x',1)$
thus $(x,1)\in X$ and we are done. Let $z$ be such an element. Let
$x\in E$. If $G= G_4(C)$, then  $(x,0)$ separates $(z,0)$ and
$(z,1)$, thus belongs to $X$. Furthermore, since $(x,1)$ separates
$(z,0)$ and $(x,1)$,  $(x,1)\in X$, proving that $X=V$. Suppose that
$G\in \{G_c(C),G_{ab}(C)\}$.  Let $\gamma\in \{c, a\}$. Since
$\gamma$ separates $(z,0)$ and $(z,1)$, it belongs to $X$. If
$\gamma=c$ then $(x, 0)$ separates $\gamma$ and $(z,0)$, hence
belongs to $X$; furthermore, since $(x,1)$ separates $\gamma$ and
$(x,0)$, it belongs to $X$, proving that $X=V$. If $\gamma=a$ then,
since it separates $(z,0)$ and $(z,1)$, it belongs to $X$. Since
$(x,0)$ separates $a$ and $b$ and $(x,1)$ separates $a$ and $(x,0)$,
the vertices $(x,0)$ and $(x,1)$ belong to $X$, proving that $X=V$.

If  $X\not = V$,  it follows from  Claim 1 and Claim 2 that $X$
contains at most one element $(x,i)$ of $E\times \{0,1\}$. Since
$\vert X\vert \geq 2$, this is impossible if $G= G_4(C)$. Suppose
that  $G\in \{G_c(C),G_{ab}(C)\}$ and let $\gamma \in X\setminus
E\times \{0,1\}$. If $\gamma \in \{c, b\}$ then  since $(x,j)$, with
$j\not = i$,  separates $\gamma$  and $(x,i)$, $(x,j)\in X$ a
contradiction. We may then suppose that $G=G_{ab}(C)$ and
$\gamma=a$. In this case, since $b$ separates $a$ and $(x,i)$, $b\in
X$, and we are lead to the previous case, which yield a
contradiction. In all cases $X=V$, hence $G$ is prime.

$(2)$ The graph $G$ is minimal prime. Let $V'\subseteq V$ such that
$\vert V'\vert= \vert V\vert$ and $G_{\restriction V'}$ is prime.
Our goal is to define an embedding from $G$ into $G'$. Set $E'(i):=
\{x\in E: (x,i)\in V'$ for $i\in \{0,1\}$. We will define $f_i:
E\rightarrow E'(i)$ ($i<2$) such that the map $F$ defined by
$F(x):=x$ for $x\in K$ and $F(x,i):=(f_i(x),i)$ is an embedding. In
order to do so, it will be enough that:
\begin{enumerate}[{(i)}]
\item $f_i(x)<f_i(y)$ for all $x<y$ and  $i<2$;
\item $f_0(x)\leq f_1(y)$ if and only if $x\leq y$.
\end{enumerate}
Suppose that $G=G_{i}(C)$ for $i\in \{0, 3,4\}$. Set $E':=E'(0)\cap
E'(1)$. Observe that the symmetric difference $E'(0)\Delta E'(1)$
has at most two elements. This yields $\vert E'\vert= \vert E\vert$.
Let $f: E\rightarrow E'$ be one to one (in this case, we do not need
to impose that $f$ is order preserving). Set  $f_i(x):=f(x)$.
Suppose that $G\in \{G_{1}(C), G_c(C), G_{ab}(C)\}$.  Notice first
that if $G=G_c(C)$ or   $G_{ab}(C)$ then $G'$ must contain $c$ or
$\{a,b\}$, hence our goal reduces to define the $f_i$'s. To do so,
we use some properties of Galois lattices (in order to avoid a
transfinite enumeration,  which requires care if $\kappa$ is
singular).  Let $R':= (E'(0), \rho, E'(1))$ where $\rho:= \{(x,y)\in
E'(0)\times E' (1): (x,0)\sim (y,1)\}$. Since $G'$ is prime, it is
point determining, hence  (iii) $R'(x)\not =R'(x)$ whenever $x\not
=x'$ and similarly (iv) $R'^{-1}(y)\not =R'^{-1}(y)$ whenever $y\not
=y'$. Next, $R'$ is Ferrers,  hence $Gal(R')$ is a chain. Due to
condition  (iii), this chain is isomorphic to a subchain of
$I(C'_{0})$; similarly  $Gal(R'^{-1})$ is isomorphic to a subchain
of $F(C'_{1})$, where $C'_{i}:=C_{\restriction E_{i}}$.  For the
simplicity of the exposition, suppose that $C$ has order type
$\kappa$. In this case,  $Gal (R')$ is a well ordered chain of order
type $\kappa' +1$ with $\kappa'\leq\kappa$.  Since  $\vert V\vert
=\kappa$,  $ \kappa'=\kappa$. Let $f_0$ be the unique order
isomorphism from $C$ onto $C'_{0}$. Define $f_1$  by choosing for
$f_1(y)$ the least element of $E'_{1}$ which is greater or equal to
$f_0(y)$.\hfill $\Box$

\begin{lemma} \label{lem:notembedding}Let $\kappa$ be  an infinite initial ordinal and $C:= (E, \leq)$ be  a chain of order
type $\kappa$ or $\kappa^*$. Then:
\begin{enumerate} [{(a)}]
\item  $K(C)\leq \overline{K(C)}$ and $G_c(C)\leq \overline {G_c(C)}$.
\item $K(C)\not \leq K(C^*)$,  $G_c(C)\not \leq G_c(C^*)$ and  $G_{ab}(C)\not \leq G_{ab} (C^*)$.
\item $G_{c}(C)$, $G_{ab}(C)$ and $\overline G_{ab}(C)$ form an antichain with respect to embeddability.
\end{enumerate}
\end{lemma}
\pre (a) We suppose that $C$ as type $\kappa$ and in fact is equal
to $\kappa$. Let  $\varphi: E\times \{0,1\}\rightarrow E\times
\{0,1\}$ defined by $\varphi(x,0):= (x,1)$ and $f(y,1):= (y+1, 0)$.
Then $\varphi$ embeds    $K(C)$ into $\overline{K(C)}$. Let
$\tilde\varphi $ be defined by $\tilde\varphi((x,i)):=\varphi((x,i))
$ and $\tilde\varphi (c):=c$. This map embeds $G_c(C)$ into
$\overline  {G_c(C)}$.

(b) $\widehat {N(K(C))}$ contains a chain made of cliques of order
type $\kappa+1$.  This is not the case for $\widehat {N(K(C^*))}$,
hence  $K(C)\not \leq K(C^*)$. The rest follows.

(c) Enough to observe that for every pair $(H_i, H_j)$ of distinct
graphs in $\{G_{c}(C), G_{ab}(C), \overline G_{ab}(C)\}$ there is a
finite graph $H_{ij}$ with $H_{ij}\leq H_i$ and $H_{ij}\not \leq
H_j$. A simple inspection shows that this can be achieved with
graphs of size at most 6.\hfill $\Box$\\

Since $G_c(C)\leq \overline {G_c(C)}$ (cf. Lemma \ref
{lem:notembedding}), we do not need to add $ \overline {G_c(C)}$ to
the  set $\mathcal M$ of graphs listed in Theorem
\ref{thm:listminimalprime graph}. To complete the proof of Theorem
\ref{thm:listminimalprime graph} we need only to prove that
$\mathcal M$ forms an antichain with respect to embeddability. We
divide it into three subsets, namely $\mathcal B$ made of those
graphs which are bipartite, $\overline {\mathcal B}$ made of the
complements of these graphs and $\mathcal R$ made of the remaining
graphs. Clearly, $\mathcal B$  is an antichain, hence $\mathcal
B\cup \overline {\mathcal B}$ is an antichain. Each member of
$\mathcal R$ is the union of an infinite independent set  and an
infinite clique (plus, possibly, an extra element), hence is
incomparable to all members of $\mathcal B\cup \overline {\mathcal
B}$. To conclude it remains to show that $\mathcal R$ is an
antichain. Since in $G_4(C)$ the members of the independent set have
degree 1, $G_4(C)$ is incomparable to the other members of $\mathcal
R$. Since the complement of a member of $\mathcal R$ embeds into an
other member, the same holds true for $\overline {G_4(C)}$. Thus, we
are left to show that the six remaining graphs   $G_c(C), G_{ab}(C),
\overline {G_{ab}(C)}$ for $C\in \{\kappa, \kappa^*\}$ form an
antichain. We may apply Lemma \ref{lem:notembedding}. But, as we
will see in the next section, these graphs are comparability graphs.
From their pictorial representation it is easy to  see that  the
twelve transitive orientations of these graphs form an antichain
(with respect to the embeddability relation between posets); in
particular these graphs form an antichain.

\subsection{Examples of minimal prime posets}

Let $C:= (E, \leq )$ be an infinite chain. The graphs $G_0(C)$ and
$G_3(C)$  are comparability graphs,  exactly as $G_0$ and $G_3$ are.
Indeed, $G_0(C)=Comp(Q_0(C))$, where $Q_0(C)$ is the set of atoms
and coatoms of  $\mathcal P(E)$ ordered by inclusion, whereas
$G_3(C)= Comp(Q_3(C))$, where $Q_3(C)$ is the set $E\times \{0,1\}$
augmented of an element $c$ and ordered  so that $(x,i)<(x,j)$ if $
i=0,j=1$ and $c<(x, 1)$ for all $x\in E$.  The  graphs   $\overline
{G_0(C)}$,  $\overline {G_3(C)}$, $G_4(C)$ and $\overline {G_4(C)}$
are not comparability graphs. The graph $G_1(C)$ and its complement
are comparability graphs.  In fact, $G_1(C)=Comp(Q_1(C))$,  where
$Q_1(C)$ is the set $E\times \{0,1\}$ ordered so that  $(x,i)<(y,j)$
if $ i=0,j=1$ and $x\leq y$,  whereas  $\overline
{G_1(C)}=Comp(P_{1}(C))$,  where $P_1(C)$ is the set $E\times
\{0,1\}$ ordered so that  $(x,i)<(y,j)$ if  $ i\geq j$ and $x> y$.

\begin{lemma}\label{invent1} If $\kappa$ is an infinite initial ordinal, the
seven posets $Q_0(\kappa)$, $Q_1(\kappa)$, $Q_1(\kappa)^*$,
$Q_3(\kappa)$, $Q_3(\kappa)^*$, $P_1(\kappa)$ and $P_1(\kappa)^*$
are minimal prime. If $\kappa=\omega$, then with  the one way
infinite fence $Q_2$ and the transitive orientations $P_2$ and
$P_2^*$ of $\overline{Q_2}$ represented Figure \ref {critique3},
they form an antichain of ten minimal prime posets.
\end{lemma}
 In order to obtain more prime posets,  we order $E\times \{0,1\}$ by
setting $(x,i)<(y,j)$ whenever $i=0$ and $x\leq y$. Let $P_c(C)$ be
the poset obtained by adding to $E\times \{0,1\}$ an extra element
$c$ in such a way that $c<(x,1)$ for all $x\in E$ and $P_{ab}(C)$
 be the poset  obtained by adding  two extra elements $a$ and $b$ to $E
\times \{0, 1\}$ in such a way that 1) $(x, 0) < a$ for all $x \in
E$ and 2) $b < a$. Let  $P_{a/b}(C)$ be the poset obtained from
$P_{ab}(C)$ by removing the comparabilities $((x,0), (x,1))$ for all
$x \in E$ and $(b,a)$ and adding the comparabilities  $(x, 1) < b$
and $(x,0)<b$ for all $x \in E$.  The posets obtained by taking $C$
to be the chain $\omega$ and then its dual $\omega^*$ are
represented in Figure \ref{critique2}. As it is easy to check, we
have:


\begin{figure}[h]
  \begin{center}
    \includegraphics[width=6in]{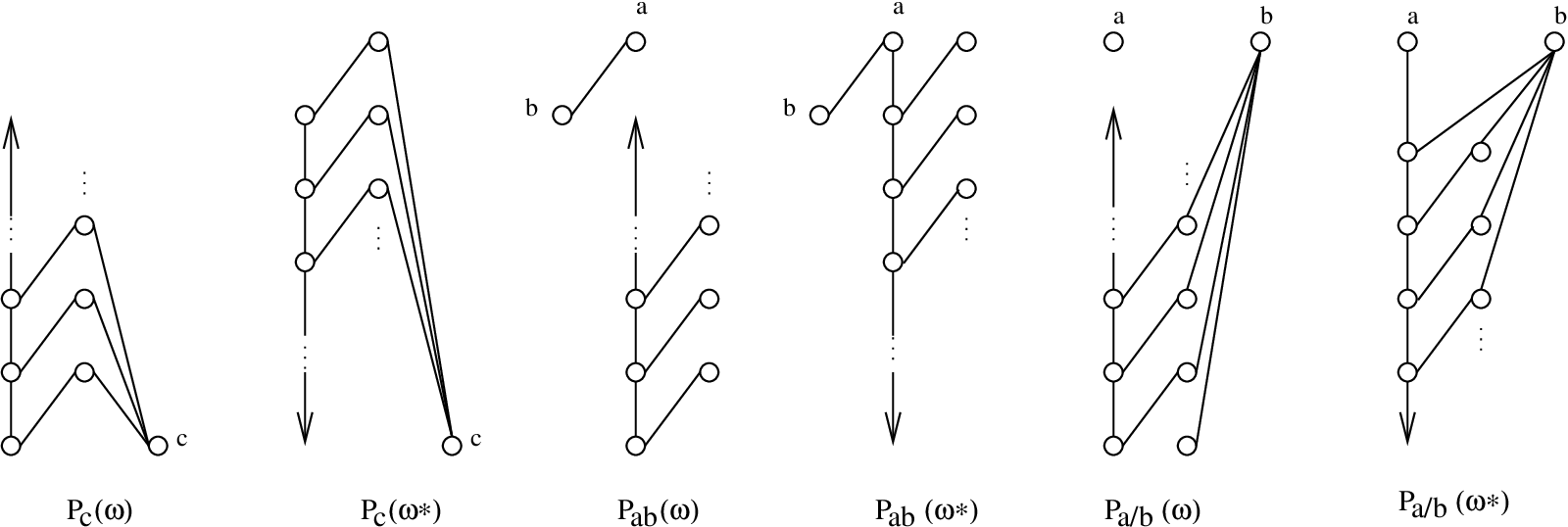}
  \end{center}

  \caption{Minimal prime posets of dimension 2.}
  \label{critique2}
\end{figure}

\begin{lemma}\label{lem:compprime} $G_c(C)=Comp(P_c(C))$, $G_{ab}(C)=Comp(P_{ab}(C))$ and
$\overline {G_{ab}(C)}$ is isomorphic to $Comp (P_{a/b}(C))$, via
the map $\varphi$ defined by $\varphi(a):=a$, $\varphi (b):=b$ and
$\varphi ((x,i)):=(x, i+1)$ (where $0+1=1$ and $1+1=0$).
\end{lemma}

Note that, from (a) of Lemma \ref{lem:notembedding}, $\overline
{G_c(C)}\leq G_c(C)$, hence $\overline {G_c(C)}$ is a comparability
graph. \begin{lemma}\label{invent2} The posets $P_c(C)$, $P_{ab}(C)$
and  $P_{a/b}(C)$, where $C\in \{\kappa, \kappa^*\}$,  are minimal
prime of dimension 2.  With  their dual, they form an antichain of
twelve minimal prime graphs.
\end{lemma}
\pre The fact that they are minimal prime follows from Lemma
\ref{lem:compprime} and Theorem \ref {thm:listminimalprime graph}.
Since the complement of their comparability graph is a comparability
graph, they have dimension 2.  The fact that they form an antichain
follows from Lemma \ref{lem:notembedding} and a careful examination
of Figure  \ref{critique2}. \hfill $\Box$

\subsection {Proof of Theorem \ref{thm:examples}\label {proofthm:examples}}
Let $\kappa$ be an infinite cardinal. Theorem
\ref{thm:listminimalprime graph} yields fourteen  minimal prime
graphs which are pairwise incomparable. The inventory made in Lemma
\ref {invent1} and Lemma \ref{invent2} of those which are
comparability graphs yields nineteen minimal prime posets. If
$\kappa=\omega$, we may add to the list of minimal prime graphs the
infinite fence and its complement and to the list of prime graphs
the infinite fence and  the two transitive orientations of its
complement.

\section{Open questions}

The countable minimal prime graphs described in Section
\ref{sample} consist of the $G_{i}$'s, for $i< 5$, and their
complements plus six graphs obtained from $K({\omega})$ and
$K({\omega^*})$ by adding one or two vertices.

\begin{question}
Does these sixteen graphs are the only countable minimal prime
graphs?
\end{question}

A preliminary question is:

\begin{question}
Does a countable prime graph embedding neither $K({\omega})$ nor
$K({\omega})$, necessarily embeds one of the graphs $G_{i}$ or
$\overline{G_{i}}$ for $i< 5$.
\end{question}

In this paper, we have described some countable minimal prime graphs
and posets. All our examples,  except one, the path,  extend to
arbitrary infinite cardinality. And so far we have obtained fourteen
minimal prime graph in each uncountable cardinality. One could try
to characterize minimal prime graphs and posets of any cardinality.

Another possible direction for future research on this subject is
the study of minimal prime relational structures. The notion of an
autonomous set for general relational structures was introduced by
Fra\"{\i}ss\'e \cite{fr} who used the term "interval" rather than
autonomous set. We can therefore define prime relational structures
in a similar way as for prime graphs and posets. But it must be
noticed that even in the case of directed graphs without circuits
the number of those which are countable and minimal prime is at
least countable. Moreover there are infinite prime directed graphs
without circuits which do not embed a countable minimal prime
directed graph. To illustrate observe that all orientations of a one
way directed graph are prime. These orientations being coded by an
infinite word on a two letter alphabet, the minimal ones are coded
by periodic words, whereas those embedding a minimal prime graph are
coded by eventually periodic words.

Still, for posets and tournaments, we ask:

\begin{questions}
\begin{enumerate}[(a)]
\item Does every infinite prime poset, respectively tournament, embeds a countable minimal prime poset, respectively tournament?
\item Are there only finitely many infinite countable minimal prime posets, respectively tournament?

\end{enumerate}
\end{questions}

Some countable minimal tournaments have been identified in
\cite{boud-pouzet}.

In the special case of posets, Theorems \ref{mainposet1} and
\ref{mainposet2} yield respectively six and four countable prime
minimal posets. Among these posets seven have dimension two. Posets
depicted in Figures \ref{critique2} and their dual yield twelve
countable minimal prime posets with dimension 2. We do not know
whether this list is complete. A preliminary question is this.

\begin{questions}

Do the nineteen posets of dimension 2 mentioned above are the only
countable minimal prime posets with dimension 2?

\end{questions}

\end{document}